\documentclass[12pt]{article}
\usepackage{amssymb}
\usepackage{amsmath}
\usepackage{latexsym}
\usepackage{pgf, tikz}
\usetikzlibrary{arrows}

\newcommand\subs\subseteq
\newcommand\emp\emptyset

\renewcommand\a{\alpha}

\renewcommand\o{\omega}
\renewcommand\phi{\varphi}

\renewcommand\hbar{{\overline{h}}}

  \newcommand\pf{\medskip\noindent{\bf Proof.}\ }

  \newcommand\dinsq{{\square \kern-0.75em \diamondsuit}}
  
  \newcommand\qed{\hfill{\vbox{\hrule\hbox{\vrule\kern3pt
                  \vbox{\kern6pt}\kern3pt\vrule}\hrule}}}

  \newcommand\EE{{\cal E}}

  \newcommand\RR{\mathbb{R}}

  \newcommand\tams{{\em Trans. Amer. Math. Soc.}}

\author{Bal\'azs Bursics \and P\'eter Komj\'ath\footnote{Partially supported by 
Hungarian National Research Grant OTKA K 131842}}

\title{A coloring of the plane without monochromatic right triangles} 

\begin{document}

  \maketitle

\begin{abstract} 
We give a full, correct proof of the following result, earlier claimed in 
\cite{kope1}.  
If the Continuum Hypothesis holds then there is a coloring of the plane 
with countably many colors, with no monocolored right triangle. 
\end{abstract}

\medskip 
One 
 of Paul Erd\H{o}s's favorite topics consisted the applications of the 
Axiom of Choice to construct paradoxical sets and colorings of the 
Euclidean spaces. 
Among a large number of other questions, he raised the following:  
is there a coloring of the plane with countably many colors, with 
no monocolored right angled triangles.  
In \cite{kope1} this was shown to be equivalent to the Continuum Hypothesis. 

Recently, the senior author observed that the proof of the positive 
direction in \cite{kope1} is incomplete. 
Here we give a full, correct proof.

We notice that later Schmerl in \cite{schmerl} gave another, more general 
proof, which, however, is less elementary.

\medskip
\noindent{\bf Notation. Definitions.} 
We use the notation and terminology of axiomatic set theory. 
The ordinals are von Neumann ordinals, $\o_1$ is the least uncountable 
ordinal.


\medskip
\noindent{\bf Theorem.} 
(CH) {\em There is $f:\RR^2\to\o$ with no monochromatic right angles.}

\pf 
Using CH, we decompose $\RR^2$ and the set $\EE$  of planar lines and circles 
into the increasing, continuous sequence of sets as 
$\RR^2=\bigcup\{H_\a:\a<\o_1\}$ and $\EE=\bigcup\{\EE_\a:\a<\o_1\}$ such 
that $H_0=\EE_0=\emp$ and \\
(1)  if $x\neq y\in H_\a$, then their connecting line and Thales circle 
are in $\EE_\a$, \\ 
(2) if $x,y,z\in H_\a$ are not collinear, then their circuit is in $\EE_\a$, \\
(3) if $e_0\neq e_1\in\EE_\a$, then $e_0\cap e_1\subs H_\a$, \\
(4) if $x\in C\in \EE_\a$, $C$ is a circle, then the antipodal point of $x$ 
in $C$ is in $H_\a$, \\ 
(5) if $L\in\EE_\a$ is a line, $x\in H_\a\cap L$, then the line $L'$ 
perpendicular to $L$ with $x\in L'$ also contained in $\EE_\a$. 

This can be done by the usual Skolem-type closing arguments. 

We are going to construct the coloring $f:\RR^2\to\o$ and the function 
$\phi:\EE\to[\o]^\o$ by transfinite recursion on $\a<\o_1$ for 
$H_{\a+1}-H_\a$ and for $\EE_{\a+1}-\EE_\a$, satisfying the following: \\ 
(6) $f|(H_{\a+1}-H_\a)$ is injective, \\ 
(7) if $x\in H_{\a+1}-H_\a$, $e\in \EE_{\a}$, $x\in e$, then $f(x)\in \phi(e)$, \\ 
(8) if $C\in\EE$ is a circle, $i\in \phi(C)$, $x,y\in C$, $f(x)=f(y)=i$, 
then $x,y$ are not antipodal, \\ 
(9) if $C\in \EE$ is a circle, $i\notin \phi(C)$, then 
$|f^{-1}(i)\cap C|\leq 2$, \\ 
(10)  if $L\in\EE$ is a line, $i\notin \phi(L)$, then 
$|f^{-1}(i)\cap L|\leq 1$, \\ 
(11) if $L,L'$ are perpendicular lines, $\{x\}=L\cap L'$, then 
$f(x)\notin \phi(L)\cap \phi(L')$. 
 
\medskip
\noindent{\bf Claim 1.} {\em There is no right triangle monocolored by $f$.} 

\pf 
Assume that $x,y,z$ form a right triangle with the right angle at $y$ and 
$i=f(x)=f(y)=f(z)$.  
Let $C$ be the circle around $x,y,z$. 
If $i\notin \phi(C)$, then we obtain a contradiction with (9). 
If $i\in \phi(C)$, then we get a contradiction with (8), as $x,z$ are 
antipodal. 
\qed 

\medskip
\noindent{\bf Claim 2.} 
         {\em If $x\in H_{\a+1}-H_\a$, then there is at most one $e\in\EE_\a$ 
         such that $x\in e$.} 

\pf 
By (3). 
\qed

\medskip
We add the following condition: \\ 
(12) if $x\in H_{\a+1}-H_\a$, $L$ is a line with $x\in L\in\EE_\a$, $L'$ is the 
line perpendicular to $L$ at $x$, $\{y\}=L'\cap H_\a$, then $f(x)\neq f(y)$.  

Notice that $L'\in \EE_{\a+1}-\EE_\a$ by (3) and (5), and $|L'\cap H_\a|\le 1$ by (1). 

Assume that $f|H_\a$ and $\phi|\EE_\a$ are already constructed, we have 
to define $f|(H_{\a+1}-H_\a)$ and $\phi|(\EE_{\a+1}-\EE_\a)$. 

Enumerate $H_{\a+1}-H_\a$ as $\{x_j:j<\o\}$. 
By recursion on $j$ define $f(x_j)$ so that 
\[
f(x_j)\geq\max\{f(x_0),\dots,f(x_{j-1})\}+2,
\]
$f(x_j)$ satisfies (7), 
and 
is different from the (possible) color disqualified by 
(12). 
Clearly (6) is satisfied, and also $\o-f[H_{\a+1}-H_\a]$ is infinite. 

Next we define $\phi$ on the circuits in $\EE_{\a+1}-\EE_\a$. 
If $C$ is such a circuit, set $A=C\cap H_\a$. 
By (2), we have $|A|\leq 2$. 

\medskip
\noindent{\bf Case 1.} 
$|A|\leq 1$ or $f$ takes distinct values on the two elements of $A$. 

In this case, set $\phi(C)=\o-f[A]$.

\medskip
\noindent{\bf Case 2.} 
$|A|=2$ and $f$ assumes the same value on the elements of $A$. 

In this case, set $\phi(C)=\o$. 

\medskip 
Finally, we define $\phi$ on the lines in $\EE_{\a+1}-\EE_\a$. 
Let $L\in\EE_{\a+1}-\EE_\a$ be a line. 
Set $B=L\cap H_\a$. 
Notice that by (1), $|B|\leq 1$.  
If $B$ is nonempty, let $y$ be its unique point.

\medskip
\noindent{\bf Case 1.} 
$B=\emp$ or $f(y)\notin f[L\cap(H_{\a+1}-H_\a)]$. 

Then let $\phi(L)=\o-f[L\cap H_{\a+1}]$. 

\medskip
\noindent{\bf Case 2.}  
$f(y)\in f[L\cap (H_{\a+1}-H_\a)]$. 

Then set $\phi(L)=(\o-f[L\cap(H_{\a+1}-H_\a])\cup \{f(y)\}$.

\medskip
\noindent{\bf Claim 3.} $\phi(e)\in[\o]^\o$ ($e\in \EE_{\a+1}-\EE_\a$). 

\pf 
If $C\in\EE_{\a+1}-\EE_\a$ is a circuit, then this is obvious in Case 2 
and in Case 1 $\phi(C)$ is $\o$ minus at most 2 elements. 

If $L\in \EE_{\a+1}-\EE_\a$ is a line, the statement follows as 
$f[H_{\a+1}-H_\a]$ is a coinfinite set.  
\qed 

\medskip
\noindent{\bf Claim 4.} $\phi$ {\em satisfies}  (8). 

\pf 
Assume that $C\in\EE_{\a+1}$ is a circle, $x,y\in C\cap H_{\a+1}$ are 
antipodal and $f(x)=f(y)=i$. 

\medskip
\noindent{\bf Case 1.} $C\in\EE_\a$.  

In  this case one of $x,y$, say $x$ must be in $H_\a$ by (6). 
Then, by (4), $y$ is also in $H_\a$, and we are finished by induction. 

\medskip
\noindent{\bf Case 2.} $C\in\EE_{\a+1}-\EE_\a$. 

If $x,y\in H_{\a+1}-H_\a$, then $f(x)\neq f(y)$ by (6). 

If $x,y\in H_\a$, then $C\in\EE_\a$ by (1), a contradiction again. 

Assume finally, that $x\in H_\a$, $y\in H_{\a+1}-H_\a$. 
If, in the definition of $\phi(C)$, Case 1 applies, then 
$i\notin \phi(C)$. 
We can therefore assume that Case 2 holds, $H_\a\cap C=\{x,z\}$ and 
$f(y)=f(z)=i$. 
Let $L$ be the connecting line of $x$ and $z$, $L'$ the connecting line of 
$z$ and $y$. 

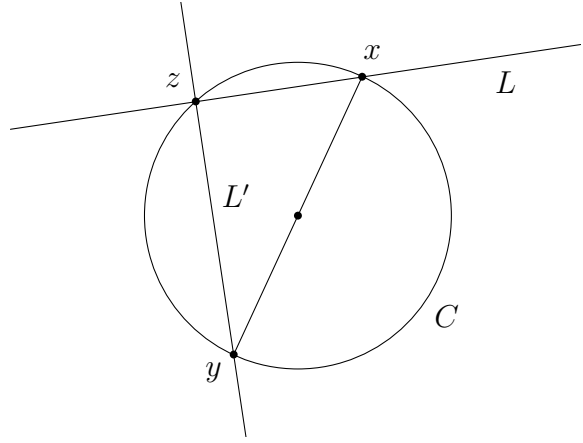
\smallskip
\begin{center}
    \begin{tikzpicture}[line cap=round,line join=round,>=triangle 45,x=0.7cm,y=0.7cm]
\clip(-3,-2.73) rectangle (7.83,5.51);
\draw(2.46,1.46) circle (2.04cm);
\draw [domain=-3:7.83] plot(\x,{(--5.11-4.81*\x)/0.72});
\draw [domain=-3:7.83] plot(\x,{(--11.21--0.47*\x)/3.16});
\draw (3.68,4.1)-- (1.24,-1.18);
\draw (3.5,4.9) node[anchor=north west] {$ x $};
\draw (0.5,-1.12) node[anchor=north west] {$ y $};
\draw (-0.27,4.38) node[anchor=north west] {$ z $};
\draw (4.85,-0.03) node[anchor=north west] {$ C $};
\draw (6,4.41) node[anchor=north west] {$ L $};
\draw (0.81,2.27) node[anchor=north west] {$ L' $};
\begin{scriptsize}
\fill [color=black] (2.46,1.46) circle (1.5pt);
\fill [color=black] (3.68,4.1) circle (1.5pt);
\fill [color=black] (1.24,-1.18) circle (1.5pt);
\fill [color=black] (0.52,3.63) circle (1.5pt);
\end{scriptsize}
\end{tikzpicture}
\end{center}

Then $L\in\EE_\a$ by (1) and then $L'\in\EE_\a$ by (5). 
Further, by (10), $i\in \phi(L)\cap \phi(L')$, which, with $f(z)=i$, 
contradicts (11).
\qed 
 
\medskip
\noindent{\bf Claim 5.} (9) {\em holds}. 

\pf 
Assume that $C\in \EE_{\a+1}$ and $i\notin \phi(C)$. 
If $C\in\EE_\a$, then $C\cap f^{-1}(i)$ does not contain element from $H_{\a+1}-H_\a$ 
by (7). 
If $C\in\EE_{\a+1}-\EE_\a$ and Case 2 holds in the definition of 
$\phi(C)$, then there is nothing to prove. 
Assume finally, that $C\in\EE_{\a+1}-\EE_\a$, and Case 1 holds. 
Then, if $i\notin\phi(C)$, then $C\cap f^{-1}(i)$ has at most one 
element in $H_\a$, at most one in $H_{\a+1}-H_\a$ by  (6), that is at most 2.  
\qed 

\medskip
\noindent{\bf Claim 6.} (10) {\em holds}. 

\pf 
Assume that $L\in\EE_{\a+1}$, $i\notin\phi(L)$. 
If $L\in \EE_\a$, then $L\cap f^{-1}(i)$ does not increase in $H_{\a+1}-H_\a$. 
If $L\in\EE_{\a+1}-\EE_\a$ and $B=L\cap H_\a$ is empty, then 
$|L\cap f^{-1}(i)|\leq 1$ by (6). 
Otherwise, $B$ is a singleton by (1), let $y$ be its unique element. 
If $f(y)\notin f[L\cap(H_{\a+1}-H_\a)]$, then again $|L\cap f^{-1}(i)|\leq 1$, otherwise we are in Case 2 of the definition of $\phi(L)$ where we 
specifically added $f(y)$ to $\phi(C)$.  
\qed 

\medskip
\noindent{\bf Claim 7.} 
         (11) {\em holds.} 

\pf 
Assume that $L,L'\in\EE_{\a+1}$ are perpendicular, $L\cap L'=\{x\}$, 
$i=f(x)\in\phi(L)\cap \phi(L')$. 

\medskip
\noindent{\bf Case 1.} 
$L,L'\in\EE_\a$. 

In this case $x\in H_\a$ by (3), so the configuration in (11) already 
appers in $H_\a$, $\EE_\a$. 

\medskip
\noindent{\bf Case 2.} 
$L\in\EE_\a$, $L'\in\EE_{\a+1}-\EE_\a$.

As $i\in\phi(L')$, by the definition of the latter there is $y\in L'\cap H_\a$, 
$f(y)=i$. 
This is exactly what is ruled out at the coloring of $x$ by (12).  

\medskip
\noindent{\bf Case 3.} 
$L,L'\in\EE_{\a+1}-\EE_\a$. 

\medskip
\noindent{\bf Subcase 3.1.} 
$x\in H_\a$. 

As $i=f(x)\in\phi(L)$, by the definition of $\phi(L)$ there is 
$y\in L\cap (H_{\a+1}-H_\a)$ with $f(y)=i$. 
Likewise, there is $z\in L'\cap (H_{\a+1}-H_\a)$ with $f(z)=i$. 
Now $y,z$ are distinct elements of $H_{\a+1}-H_\a$ and $f(y)=f(z)$, 
contradicting (6). 

 \medskip
\noindent{\bf Subcase 3.2.} 
$x\in H_{\a+1}-H_\a$. 
 
As $i=f(x)\in\phi(L)$, there is $y\in L\cap H_\a$, $f(y)=i$. 
Similarly, there is $z\in L'\cap H_\a$, $f(z)=i$. 

Let $C$ be the circuit containing $x,y,z$. 

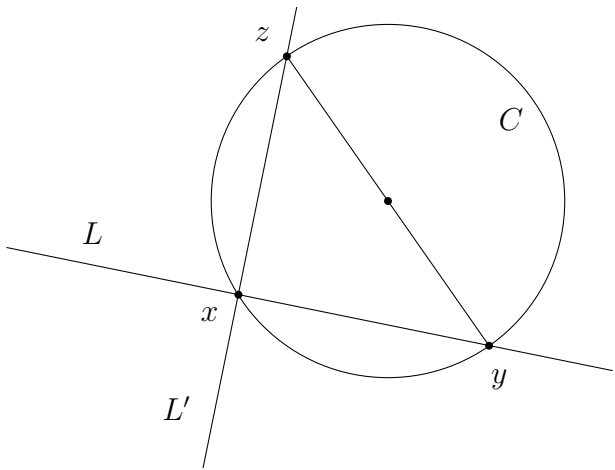
\begin{center}
    \begin{tikzpicture}[line cap=round,line join=round,>=triangle 45,x=0.7cm,y=0.7cm]
\clip(-5.41,-3.17) rectangle (6.09,5.58);
\draw(1.83,1.9) circle (2.35cm);
\draw [domain=-5.41:6.09] plot(\x,{(-0.4-0.97*\x)/4.77});
\draw [domain=-5.41:6.09] plot(\x,{(--4.7--4.53*\x)/0.92});
\draw (-0.09,4.65)-- (3.75,-0.85);
\draw (-1.93,0.08) node[anchor=north west] {$ x $};
\draw (3.58,-1.1) node[anchor=north west] {$ y $};
\draw (-0.91,5.4) node[anchor=north west] {$ z $};
\draw (3.73,3.87) node[anchor=north west] {$ C $};
\draw (-4.19,1.71) node[anchor=north west] {$ L $};
\draw (-2.66,-1.61) node[anchor=north west] {$ L' $};
\begin{scriptsize}
\fill [color=black] (1.83,1.9) circle (1.5pt);
\fill [color=black] (-0.09,4.65) circle (1.5pt);
\fill [color=black] (3.75,-0.85) circle (1.5pt);
\fill [color=black] (-1.01,0.12) circle (1.5pt);
\end{scriptsize}
\end{tikzpicture}
\end{center}

As $L,L'$ are perpendicular, $y$ and $z$ are antipodal in $C$. 
As in $C$ there are 3 points of color $i$, $i\in \phi(C)$. 
But this contradicts to the antipodality of $y$ and $z$. 
\qed

\medskip
The proof of the theorem is concluded. 
\qed

\bigskip
\hbox{
\vtop{\hbox{Bal\'azs Bursics}
      \hbox{Institute of Mathematics}
      \hbox{E\"otv\"os University}
      \hbox{Budapest, P\'azm\'any P. s. 1/C} 
      \hbox{1117, Hungary}
      \hbox{e-mail:{\tt\ bursicsb@gmail.com}}}}

\bigskip
\hbox{
\vtop{\hbox{P\'eter Komj\'{a}th}
      \hbox{Institute of Mathematics}
      \hbox{E\"otv\"os University}
      \hbox{Budapest, P\'azm\'any P. s. 1/C} 
      \hbox{1117, Hungary}
      \hbox{e-mail:{\tt\ peter.komjath@gmail.com}}}}

\end{document}